\documentclass[12pt,thmsa]{amsart}

\def\Cal{\mathcal}

\def\P{{\Cal P}}

\def\S{{\Cal S}}
\def\F{{\Cal F}}

\def\I{{\Cal I}}
\def\L{{\Cal L}}

\def\W{{\Cal W}}

\def\tr{{\hbox{\rm tr}}}

\def\Ma{\frM_{n,m}}
\def\Mt{\frM_{n-k,m}}
\def\Mkm{\frM_{k,m}}

\def\W{\mathcal{W}}

\def\f0{f_0}
\def\Fc0{\varphi_0}
\def\rn{\bbr^n}

\def\rnk{\bbr^{n-k}}

\def\I_k {I_{-}^{k/2}}
\def\I+k {I_{+}^{k/2}}

\def\vnk{V_{n,n-k}}
\def\vnm{V_{n,m}}
\def\cd{\vnk\times\frM_{n-k,m}}
\def\Gr{G(n,k,m)}
\def\sigk{\sig_{n,n-k}}

\def\bbr{{\Bbb R}}

\def\bbn{{\Bbb N}}

\def\bbc{{\Bbb C}}

\def\rank{{\hbox{\rm rank}}}
\def\diag{{\hbox{\rm diag}}}

\def\mod{{\hbox{\rm mod}}}

\def\tr{{\hbox{\rm tr}}}

\def\det{{\hbox{\rm det}}}

\def\min{{\hbox{\rm min}}}

\def\fc{\varphi(\xi,t)}

\def\lp{L^p(\frM_{n,m})}

\def\part{\partial}
\def\intl{\int\limits}
\def\b{\beta}

\def\Gam{\Gamma}

\def\a{\alpha}

\def\Del{\Delta}

\def\vp{\varphi}

\def\g{\gamma}
\def\gam{\gamma}

\def\sig{\sigma}
\def\lam{\lambda}

\def\t{\tau}
\def\eq{\xi 'x=t}
\def\eqv{\bar\xi '\bar x=\bar t}
\def\Rnm{\bbr^{nm}}

\def\Rnkm{\bbr^{(n-k)m}}
\def\rf{\hat f (\xi,t)}
\def\df{\check \varphi (x)}

\font\frak=eufm10

\def\fr#1{\hbox{\frak #1}}

\def\frL{\fr{L}}        
\def\frM{\fr{M}}

\def\gma{\Gamma_m(\a)}

\def\gk{\Gamma_k}
\def\lin{{\hbox{\rm lin}}}

\def\det{{\hbox{\rm det}}}

\def\min{{\hbox{\rm min}}}

\def\p{\P_m}
\def\gm{\Gamma_m}
\def\tr{{\hbox{\rm tr}}}
\def\part{\partial}
\def\intl{\int\limits}
\def\b{\beta}

\def\Gam{\Gamma}

\def\a{\alpha}

\newtheorem{theorem}{Theorem}[section]
\newtheorem{lemma}[theorem]{Lemma}

\theoremstyle{definition}
\newtheorem{definition}[theorem]{Definition}

\theoremstyle{corollary}
\newtheorem{corollary}[theorem]{Corollary}
\theoremstyle{remark}
\newtheorem{remark}[theorem]{Remark}
\usepackage{amssymb}

 \numberwithin{equation}{section}

\newcommand{\be}{\begin{equation}}
\newcommand{\ee}{\end{equation}}

\newcommand{\bea}{\begin{eqnarray}}
\newcommand{\eea}{\end{eqnarray}}
\newcommand{\Bea}{\begin{eqnarray*}}
\newcommand{\Eea}{\end{eqnarray*}}

\begin{document}

\title [ The Fuglede formula ]{ An analogue of the Fuglede formula  in integral geometry on matrix spaces}

\author [E. Ournycheva and B. Rubin]{Elena Ournycheva and Boris Rubin}

\address{Institute of Mathematics,
Hebrew University, Jerusalem 91904,  ISRAEL}
\email{boris@math.huji.ac.il, ournyce@math.huji.ac.il}

\thanks{ The work  was supported in part by
the Edmund Landau Center for Research in Mathematical Analysis and
Related Areas, sponsored by the Minerva Foundation (Germany).}





\begin{abstract}
The  well known formula  of B. Fuglede expresses the mean value of
the
 Radon k-plane transform on $\bbr^n$ as a Riesz
 potential. We extend this formula to the space  of  $n \times m$
 real
 matrices and  show that the corresponding matrix  k-plane transform
$f \to \hat f$ is injective if and only if $n-k \ge m$. Different
 inversion formulas for this transform are obtained.
 We assume that $f \in L^p$ or $f$ is a continuous function satisfying certain ``minimal'' conditions at
 infinity.
\end{abstract}

\maketitle

\section{Introduction}

\setcounter{equation}{0}

In  1958,  B. Fuglede [Fu] proved the following remarkable formula
\be\label{fu1} c(\hat f)^{\vee} (x)=(I^k f)(x), \qquad x \in
\bbr^n, \ee where  $\hat f \equiv\hat f (\t)$ is the integral of
$f(x)$ over a $k$-dimensional plane  $\t, \; 0<k<n$; $(\hat
f)^{\vee} (x)$
 is the mean value of $\hat f (\t)$ over all $k$-planes
through $x$,  $(I^k f)(x)$  denotes  the Riesz potential of order
$k$, and $c=c(n,k)$ is a constant. For $k=n-1$, when $\t$ is a
hyperplane, this  formula was implicitly exhibited by J. Radon [R]
who was indebted to W. Blaschke for this idea. A consequence of
(\ref{fu1}) is the inversion formula \be\label{inv} f=c
(-\Del)^{k/2}(\hat f)^{\vee} \ee in which $\Del$ denotes the
Laplace operator.

Our aim is to extend (\ref{fu1}) and (\ref{inv}) to the case when
$x$ is an $n\times m$ real matrix. We note that
 Riesz potentials of functions of matrix argument and their
generalizations  arise in different contexts in harmonic analysis,
integral geometry, and PDE; see \cite{Far}, [FK], [Ge], [Kh],
[Ra], [St1], [Sh1]-[Sh3]. A  systematic study of Radon transforms
on matrix spaces  was initiated
 by E.E. Petrov [Pe1] and continued in $[\check{\rm C}]$, [Gr], [Pe2]-[Pe4], [Sh1], [Sh2].
 These
 publications  traditionally  deal with Radon transforms $f \to \hat f$ of  $C^\infty$ rapidly decreasing functions, and employ  decomposition in plane waves
 in order to recover $f$ from  $\hat f$.

We suggest another approach which
is  based on  a matrix analogue of (\ref{fu1}) (see Theorem 5.4) and allows to handle arbitrary
 continuous or locally integrable functions $f$  subject to  mild restrictions at
infinity. These restrictions are  minimal in a certain sense.
We show that inequality $n-k \ge m$ is necessary and sufficient
for injectivity of the matrix $k$-plane transform, and obtain two
inversion formulas (in terms of the Fourier transform and in the
form (\ref{inv})).

 Our key motivation is the following. Converting integral-geometrical
 entities  into  their matrix counterparts  has specific ``higher rank"
 features and sheds new light on
classical
 problems in multidimensional spaces (just note  the space
 $\bbr^N$ with $N=nm$ can be treated as a collection of  $n \times m$
 matrices).

{\bf Acknowledgements.} We are grateful to Prof. E.E. Petrov and
Dr. S.P.  Khekalo  for  the papers they sent us, and useful
discussions.

\section{Preliminaries}

\setcounter{equation}{0}

 We establish some notation and recall basic facts. Let $\frM_{n,m}$ be the
space of real matrices $x=(x_{i,j})$ having $n$ rows and $m$
 columns. We identify $\frM_{n,m}$
 with the real Euclidean space $\bbr^{nm}$ and set $dx=\prod^{n}_{i=1}\prod^{m}_{j=1}
 dx_{i,j}$. In the following
  $x'$ denotes the transpose of  $x$, $I_m$ is the identity $m \times m$
  matrix,  $0$ stands for zero entries. Given square matrices $a, b, q, s$, etc.,
  we denote by $|a|, |b|, |q|, |s|$ their determinants.   $GL(m,\bbr)$ is the group of
 real non-singular $m \times m$ matrices; $SO(n)$ is the group of
orthogonal $n \times n$ real matrices of determinant one, endowed
with the normalized invariant measure. We denote by  $\p$ the cone
of positive definite symmetric matrices $r=(r_{i,j})$;  $\tr (r)$
is the trace of $r, \; dr=\prod_{i \le j} dr_{i,j}$. The
 Lebesgue space $\lp$
 and the Schwartz space  $\S(\Ma)$ are
 identified with  respective spaces on $\bbr^{nm}$.

 \begin{lemma}\label{12.2} {\rm (see, e.g.,
 [Mu, pp. 57--59])}.\hskip10truecm

\noindent
 {\rm ($i$)} \ If $ \; x=ayb$, where $y\in\Ma, \; a\in  GL(n,\bbr), \; b \in  GL(m,\bbr)$, then
 $dx=|a|^m |b|^ndy.$\\
 {\rm ($ii$)} \ If $ \; r=q'sq$, where $s\in \P_m, \; q\in  GL(m,\bbr)$,
  then $dr=|q|^{m+1}ds.$ \\
  {\rm ($iii$)} \ If $ \; r=s^{-1}$, where $s\in \p$,   then $r\in
  \p$,
  and $dr=|s|^{-m-1}ds.$

 \end{lemma}

 The   Siegel gamma function
 associated to the cone $\p$ is defined by
\be\label{2.4}
 \gm (\a)=\intl_{\p} \exp(-\tr (r)) |r|^{\a -d} dr, \qquad d=(m+1)/2.\ee
This integral  converges absolutely
 if and only if $Re \, \a>d-1$, and represents a product of
 ordinary $\Gamma$-functions:
\be\label{2.5}
 \gm (\a)=\pi^{m(m-1)/4} \Gam (\a) \Gam (\a- \frac
{1}{2}) \ldots \Gam (\a- \frac {m-1}{2}). \ee If $1  \le k<m, \;
k \in \bbn, $ then  \be\label{2.5.1} \gma=\pi^{k(m-k)/2}
\Gamma_{k}(\a)\Gamma_{m-k}(\a-k/2), \ee \be\label{2.5.2}
\frac{\gm(\a)}{\gm(\a+k/2)}=\frac{\gk(\a+(k-m)/2)}{\gk(\a+k/2)}.\ee

For $n\geq m$, let $\vnm= \{v \in \frM_{n,m}: v'v=I_m \}$
 be the Stiefel manifold of orthonormal $m$-frames in $\bbr^n$;
  $V_{n,n}=O(n)$ is the orthogonal group in $\bbr^n$.
 We fix invariant measure $dv$ on
 $\vnm$ [Mu, p. 70]
  normalized by \be\label{2.16} \sigma_{n,m}
 \equiv \intl_{\vnm} dv = \frac {2^m \pi^{nm/2}} {\gm
 (n/2)}\;. \ee

\begin{lemma}\label{l2.3} {\rm (polar decomposition; see, e.g.,  [Mu, pp. 66,
591], [Ma]).}\hfil

\noindent Let $x \in \frM_{n,m}, \; n \ge m$. If  $\rank (x)= m$,
then
\[ x=vr^{1/2}, \qquad v \in \vnm,   \qquad r=x'x \in\p,\] and
$dx=2^{-m} |r|^{(n-m-1)/2} dr dv$.
\end{lemma}

\section{Riesz potentials}

\setcounter{equation}{0}

 For $x \in \Ma$, let
 $|x|_m =\det (x'x)^{1/2}$. The Riesz potential of order
$\a \in \bbc$ is defined by \be\label{rie} (I^\a
f)(x)=\frac{1}{\gam_{n,m} (\a)} \intl_{\Ma} f(x-y) |y|^{\a-n}_m
dy, \ee
\[ \gam_{n,m} (\a)=\frac{2^{\a m} \, \pi^{nm/2}\, \Gam_m (\a/2)}{\Gam_m
((n-\a)/2)},   \] provided this
 expression is finite. Note that for $m \ge 2$, the sets of poles of $\Gam_m (\a/2)$ and $\Gam_m
((n-\a)/2)$ are $\{m-1, \, m-2, \ldots \}$, and $\{n-m+1, \,
n-m+2, \ldots \}$ respectively. These sets overlap if and only if
$2m \ge n+2$ (keep this inequality in mind!).
 For $m=1$, the set of poles of $\Gam_m((n-\a)/2)$ is
  $\{n, \,
n+2, \, n+4, \ldots  \}$. Thus if
\[\a = \left \{
\begin{array} {ll}\displaystyle{n-m+1, \,  n-m+2, \ldots   } &
 \mbox{ for    $m \ge 2$,} \\
  \displaystyle{n, \,  n+2, \ldots }
 &\mbox { for   $m=1$}, \\
\end{array}
\right.\]
 then the  coefficient $1/\gam_{n,m} (\a)$ is
 infinite. In the following we exclude these values of $\a$ and focus on the
 case $m \ge 2$.

 The Riesz distribution corresponding to (\ref{rie}) is defined by
 \be\label{dis} (h_\a, f)= a.c. \intl_{\Ma} h_\a (x) f(x) dx,
 \qquad  h_\a (x)=\frac{|x|^{\a-n}_m}{\gam_{n,m} (\a)}, \ee
 where ``$a.c.$" abbreviates analytic continuation in the
 $\a$-variable, and $f \in \S(\Ma)$. The following lemma resumes basic properties of
 $h_\a$.

\begin{lemma}\label{rpr}{} Let $m \ge 2$. \hfill

\noindent
{\rm (i)}  The integral in (\ref{dis}) absolutely converges if and
only if $Re \, \a>m-1$.

\noindent
{\rm  (ii)} The function $\a \to  h
_\a$ extends to all $\a \in
\bbc$ as a meromorphic $\S'$-distribution with the only poles at
the points $\a=n-m+1, n-m+2, \ldots$. The order of these poles is
the same as in $\Gam_m ((n-\a)/2)$.

\noindent {\rm (iii)} $h_\a$ is a positive measure if and only if
 $\a$ belongs to the set \be\label{wa} \W \! = \! \{0,
1, 2, \ldots,  k_0 \} \cup \{\a  : Re \, \a \!
> \! m \! - \! 1; \; \a  \! \neq \!  n \! - \! m (\mod \,  1) \}, \ee
\[ k_0=\min(m-1, n-m).\]

\noindent {\rm (iv)} If $\a \neq n-m+1, n-m+2, \ldots $, and
\be\label{ft} (\F f)(y)=\intl_{\Ma}  \exp(\tr(iy'x)) f (x) dx \ee
is the Fourier transform of $f \in \S(\Ma)$, then
 \be\label{fou}
(h_\a, f)=(2\pi)^{-nm}(|y|_m^{-\a},(\F f)(y)) \ee   in the sense of analytic
continuation. In particular, \be\label{nul}
  (h_0,f)= f(0). \ee
\end{lemma}

These statements can be found in [Kh]; see also [Sh1], [Sh3],
[FK]. The set (\ref{wa}) is an analog of  the  Wallach set  in [FK], p.
137. We adopt this name  for our case too. The formula (\ref{fou})
 is known as a functional equation for the  zeta function $\a \to (h_\a, f)$.
Different proofs of  (\ref{fou}) (usually in a more general
set-up) can be found in [Ge], Chapter IV; [Ra], Prop. II-9; [Sh1],
[Far], [FK], Theorem XVI.4.3; see also [St1].

In the case $2m<n+2$, when poles of $ \Gam_m (\a/2)$ and $\Gam_m
((n-\a)/2)$ do not overlap, a simple proof of (\ref{fou}) can be
given following slight modification of the argument from [St2],
Chapter III, Sec. 3.4. Let $e_s (x)=\exp (\tr(-xsx')/4\pi), \; s
\in \p$. By the Plancherel formula, \be\label{pl}
|s|^{-n/2}\intl_{\Ma} (\F f)(y) \exp (\tr (-\pi y
s^{-1}y'))dy=\intl_{\Ma}  f(x) e_s (x) dx. \ee Multiplying
(\ref{pl}) by $|s|^{(n-\a)/2 -d}, \; d= (m+1)/2$, and integrating
in $s \in \p$, after changing the order of integration we obtain
\[ \intl_{\Ma} (\F f)(y) a (y) dy=\intl_{\Ma}  f(x) b(x)
dx, \] where \bea a (y)&=&\intl_{\p} |s|^{-\a/2 -d} \exp (\tr
(-\pi y s^{-1}y')) ds \qquad (s=t^{-1})\nonumber \\ &=& \intl_{\p}
|t|^{\a/2 -d} \exp (\tr (-\pi ty'y)) dt=\Gam_m (\a/2)\pi ^{-\a
m/2}|y|_m^{-\a}   \nonumber \eea if $Re \,\a > m-1$, and
\[ b(x)=\intl_{\p} |s|^{(n-\a)/2 -d} e_s (x) ds= \frac{(4 \pi)^{m(n-\a)/2} \,
\Gam_m ((n-\a)/2)}{ |x|^{n-\a}_m} \] if $Re \,\a<n-m+1$. A simple
computation of these integrals is performed using (\ref {2.4}) and
Lemma \ref{12.2}. Thus (\ref{fou}) follows if the set $ m-1< Re
\,\a<n-m+1$ is not vacuous, i.e., $2m<n+2$. We note that if $2m
\ge n+2$ then distributions in both sides of (\ref{fou}) are not
regular simultaneously.

Explicit representations for $h_\a$ in the discrete part of the
Wallah set (\ref{wa}) play a vital role in our consideration.
\begin{lemma}\label{ldes}  Let $f \in \S(\Ma), \; k_0=\min(m-1, n-m), \; m \ge 2$. Then for $k=1, 2, \ldots,
k_0$, \be\label{des}(h_k, f)=c\intl_{\frM_{k,m}}du \intl_{SO(n)}f
\left (\gam  \left[\begin{array} {c} u \\ 0
\end{array} \right]  \right ) \, d\gam, \ee
\be\label{con1} c=2^{-km} \,\pi^{-km/2} \, \Gam_m\Big (
\frac{n-k}{2}\Big ) / \Gam_m\Big ( \frac{n}{2}\Big ). \ee
\end{lemma}
\begin{proof} We split $x \in \Ma$ in two blocks $x=[y; b]$ where $
y \in \frM_{n,k}$ and $ b \in \frM_{n, m-k}$. Then for $Re \,
\a>m-1,$
\[ (h_\a, f)=\frac{1}{\gam_{n,m} (\a)} \intl_{\frM_{n,k}} dy
\intl_{\frM_{n,m-k}}f([y; b])\left |\begin{array}{ll}  y'y & y'b \\
b'y& b'b
\end{array}\right|^{(\a-n)/2} db \]
where $\left |\begin{array}{ll}  * & * \\
*&*
\end{array}\right|$ denotes the determinant of the respective
matrix $\left [\begin{array}{ll}  * & * \\
*&*
\end{array}\right]$. By passing to polar coordinates (see Lemma \ref{l2.3})
$y=vr^{1/2}, \; v \in V_{n,k}$, $ r \in \P_k$, we have \bea (h_\a,
f)&=&\frac{2^{-k}}{\gam_{n,m} (\a)} \intl_{V_{n,k}} dv
\intl_{\P_k} |r|^{(n-k-1)/2} dr \nonumber \\ &\times&
\intl_{\frM_{n,m-k}}f([vr^{1/2}; b])\left |\begin{array}{ll}  r & r^{1/2}v'b \\
b'vr^{1/2}& b'b
\end{array}\right|^{(\a-n)/2} db \nonumber \\
&=&\frac{2^{-k} \,\sig_{n,k} }{\gam_{n,m} (\a)} \intl_{SO(n)}d\gam
\intl_{\P_k} |r|^{(n-k-1)/2} dr  \nonumber \\ &\times&
\intl_{\frM_{n,m-k}}f_\gam ([\lam_0 r^{1/2}; b])\left |\begin{array}{ll}  r & r^{1/2}\lam'_0 b \\
b'\lam_0 r^{1/2}& b'b
\end{array}\right|^{(\a-n)/2} db. \nonumber \eea
Here \[ \lam_0= \left[\begin{array} {c} I_k \\ 0
\end{array} \right]  \in V_{n,k}, \qquad f_\gam (x)=f (\gam x).\]
We write \[ b=\left[\begin{array} {c} b_1 \\ b_2
\end{array} \right], \qquad b_1 \in \frM_{k,m-k}, \qquad b_2 \in
\frM_{n-k,m-k}.\] Since $ \lam'_0 b=b_1$, then \bea (h_\a,
f)&=&\frac{2^{-k} \,\sig_{n,k} }{\gam_{n,m} (\a)}
\intl_{SO(n)}d\gam \intl_{\P_k} |r|^{(n-k-1)/2} dr
\intl_{\frM_{k,m-k}} db_1 \nonumber \\
&\times& \intl_{\frM_{n-k,m-k}}f_\gam \left ( \left [\begin{array}{ll}   r^{1/2}&  b_1 \\
0 & b_2
\end{array}\right]\right )\, \left |\begin{array}{ll}  r & r^{1/2} b_1 \\
b'_1 r^{1/2}& b'_1 b_1 +b'_2 b_2
\end{array}\right|^{(\a-n)/2} db_2. \nonumber \eea
Note that \[ \left [\begin{array}{ll}  r & r^{1/2} b_1 \\
b'_1 r^{1/2}& b'_1 b_1 +b'_2 b_2
\end{array}\right]=\left [\begin{array}{ll}  r & 0 \\
b'_1 r^{1/2}& I_{m-k}
\end{array}\right] \left [\begin{array}{ll}  I_{k}& r^{-1/2} b_1 \\
0 & b'_2 b_2
\end{array}\right], \] and  \[ \left |\begin{array}{ll}  r & r^{1/2} b_1 \\
b'_1 r^{1/2}& b'_1 b_1 +b'_2 b_2
\end{array}\right|=\det (r) \det (b'_2 b_2),\] see, e.g., [Mu], p.
577. Therefore, \be \label{an} (h_\a, f)=c_\a \intl_{SO(n)}d\gam
\intl_{\P_k} |r|^{(\a-k-1)/2} dr \intl_{\frM_{k,m-k}} \psi_{\a -k}
(\gam, r, b_1) db_1,\ee where
\[ c_\a=\frac{2^{-k} \, \sig_{n,k} \, \gam_{n-k,m-k} (\a -k)}{\gam_{n,m}
(\a)}, \] \bea  \psi_{\a -k} (\gam, r, b_1)&=&\frac{1}{\gam_{n-k,m-k}
(\a -k)}
\intl_{\frM_{n-k,m-k}}f_\gam \left ( \left [\begin{array}{ll}   r^{1/2}&  b_1 \\
0 & b_2
\end{array}\right]\right ) \nonumber \\ &\times&  |b'_2 b_2|^{\frac{(\a-k)-(n-k)}{2}} db_2. \nonumber
\eea The last expression is the Riesz distribution of order $\a-k$
in the $b_2$-variable. Owing to (\ref{nul}), analytic continuation
of (\ref {an}) at $\a=k$ reads
\[(h_\a, f)=c_k \intl_{SO(n)}d\gam
\intl_{\P_k} |r|^{-1/2} dr \intl_{\frM_{k,m-k}} f_\gam \left ( \left [\begin{array}{ll}   r^{1/2}&  b_1 \\
0 & 0 \end{array}\right]\right ) db_1,\] $c_k=\lim\limits_{\a \to
k}c_\a$. To transform this expression, we replace $\gam$ by $\gam \left [\begin{array}{ll}   \b& 0 \\
0 & I_{n-k} \end{array}\right]$, $ \b \in SO(k)$ and integrate in $\b$. This gives \bea
(h_\a, f)&=&c_k  \! \intl_{\P_k} \! |r|^{-1/2} dr \! \intl_{SO(k)}
\! d\b  \! \intl_{\frM_{k,m-k}} \! db_1 \!
\intl_{SO(n)}  \! f_\gam \left ( \left [\begin{array}{ll}   \b r^{1/2}&  \b b_1 \\
0 & 0
\end{array}\right]\right ) d\gam \nonumber \\
&& \text{\rm (set $\zeta=\b b_1, \quad \eta=b |r|^{1/2}$ and use
Lemma \ref{l2.3})}
\nonumber \\
&=&\frac{2^{k} \, c_k}{\sig_{k,k}}\intl_{\frM_{k,k}} d\eta
\intl_{\frM_{k,m-k}} d\zeta\intl_{SO(n)} f_\gam \left ( \left [\begin{array}{ll}   \eta &  \zeta \\
0 & 0
\end{array}\right]\right ) d\gam \nonumber \\
&=&c\intl_{\frM_{k,m}}du \intl_{SO(n)}f \left (\gam
\left[\begin{array} {c} u \\ 0
\end{array} \right]  \right ) \, d\gam, \nonumber
\eea
\[ c=\frac{\sig_{n,k}}{\sig_{k,k}}  \lim\limits_{\a \to
k}  \frac{\gam_{n-k,m-k} (\a -k)}{\gam_{n,m} (\a)}=2^{-km}
\,\pi^{-km/2} \, \Gam_m\Big ( \frac{n-k}{2}\Big ) / \Gam_m\Big (
\frac{n}{2}\Big ).\] (here we used formulae (\ref{2.5.1}) and
(\ref{2.5.2})).
\end{proof}

\begin {definition}\label{cor} According to Lemma \ref{ldes} and
(\ref{nul}), we can redefine the Riesz potential $I^\a f$ for any
locally integrable function $f$  as

\be\label{any} (I^\a f)(x) = \left \{
\begin{array} {ll}\displaystyle{\frac{1}{\gam_{n,m} (\a)} \intl_{\Ma} f(x-y) |y|^{\a-n}_m
dy} \\
 \mbox{ if   $Re \, \a > m-1; \quad \a \neq n-m+1, n-m+2, \ldots $}, \\
{} \\
  \displaystyle{c\intl_{\frM_{k,m}}du
\intl_{SO(n)}f \left (x-\gam  \left[\begin{array} {c} u \\ 0
\end{array} \right]  \right ) \, d\gam} \quad \mbox { if   $\a=1, \ldots, k_0$}, \\
\end{array}
\right.
 \ee
 Here $m \ge 2, \; c$ is the constant (\ref{con1}), $k_0=\min (m-1, n-m)$. It is assumed that $f$ is
good enough, so that the corresponding integrals are absolutely
convergent.
\end{definition}

\noindent {\bf Conjecture 3.4.} We state the following hypotheses.

\noindent (a) For $f \in L^p (\Ma)$, the integrals in (\ref{any})
absolutely converge if and only if \be \label{conj} 1\leq
p<(n+m-1)/(Re \, \a+m-1). \ee (b) If $f$ is a continuous function
satisfying $ f(x)=O(|I_m+x'x|^{-\lam/2})$, then the  absolute
convergence of (\ref{any}) holds if and only if $ \lam>Re
\,\a+m-1$.

 This conjecture is true if
$\a$ is a positive integer; see  the next section.

\noindent {\bf Remark 3.5.}  Another formula for $h_k$ obtained in
[Sh1] and  [Kh] reads \be\label{rep2} (h_k, f) \! = \! c_1
\intl_{\frM_{n,k}} \! \frac{dy}{|y|_k^{n-m}}\intl_{\frM_{k,m-k}}
\!  \!  f([y; yz]) dz, \quad k \! = \! 1,2, \ldots, k_0.\ee
\be\label{con2} c_1= 2^{-km}\pi^{k(k-n-m)/2} \,
\Gam_k\left(\frac{n-m}{2}\right)/\Gam_k\left(\frac{k}{2}\right).
\ee
 It can be derived from (\ref{des}). Indeed,

\bea  (h_k, f)&=&c\intl_{\frM_{k,m}}du \intl_{SO(n)}f \left (\gam
\left[\begin{array} {c} u \\ 0
\end{array} \right]  \right ) \, d\gam \nonumber \\
&=&c\intl_{\frM_{k,m}}du \intl_{SO(n)}f(\gam\lam_0 u)d\gam \qquad
\left ( \lam_0= \left[\begin{array} {c} I_k \\ 0
\end{array} \right]  \in V_{n,k} \right )\nonumber \\
&=&\frac{c}{\sig_{n,k}}\intl_{\frM_{k,m}}du \intl_{V_{n,k}}f(v
u)dv.\nonumber \eea Now we represent $u$ in the block form
$u=[\eta; \zeta], \; \eta \in \frM_{k,k}, \; \zeta \in
\frM_{k,m-k}$, and change the variable $\zeta=\eta z$. This gives
\[ (h_k, f)=\frac{c}{\sig_{n,k}}\intl_{\frM_{k,k}}
|\eta|^{m-k}d\eta\intl_{\frM_{k,m-k}}dz \intl_{V_{n,k}}f(v[\eta;
\eta z])dv.\] Using Lemma \ref{l2.3} repeatedly, and changing
variables, we obtain \bea  (h_k, f)&=&\frac{c \, \sig_{k,k} \,
2^{-k}}{\sig_{n,k}}\intl_{\P_k} |r|^{(m-k-1)/2}
dr\intl_{\frM_{k,m-k}}dz \intl_{V_{n,k}}f(v[r^{1/2}, r^{1/2}z])dv\nonumber \\
&=&c_1\intl_{\frM_{n,k}}
\frac{dy}{|y|_k^{n-m}}\intl_{\frM_{k,m-k}} f([y; yz]) dz\nonumber
\eea
where by  (\ref{con1}), (\ref{2.16}) and (\ref{2.5.2}),
\[ c_1=\frac{c \, \sig_{k,k} }{\sig_{n,k}}=
 2^{-km}\pi^{k(k-n-m)/2} \,
\Gam_k\left(\frac{n-m}{2}\right)/\Gam_k\left(\frac{k}{2}\right).\]

\section{Radon transforms}\label{s4}

\setcounter{equation}{0}

\subsection{Matrix planes}\label{s4.1}
Let $k,n$, and $m$ be positive integers, $0<k<n$, $\vnk$ be the
Stiefel manifold of orthonormal $(n-k)$-frames in $\bbr^n$. For
$\; \xi\in\vnk$, $t\in\Mt$, the linear manifold \be\label{plane}
\tau= \tau(\xi,t)=\{x\in\Ma:\eq\} \ee
 will be called {\it a matrix $k$-plane} in $\Ma$. We denote by  $\Gr$  a variety of all such
 planes. The parameterization $\tau=
\tau(\xi,t)$ by the points $(\xi,t)$ of the ``matrix cylinder" $\cd$
is not one-to-one because for any orthogonal transformation $\theta \in
O(n-k)$, the pairs $(\xi,t)$ and $(\xi\theta ', \theta t)$ define the same
plane. We identify  functions $\varphi(\t)$ on $\Gr$ with
functions $\fc$ on $\cd$ satisfying $\varphi(\xi\theta ',\theta t)=\fc$
for all $\theta\in O(n-k)$, and supply  $\Gr$ with the measure $d\tau$
so that \be \intl_{\Gr} \varphi(\t) \, d\tau=\intl_{\cd}
 \fc \, d\xi dt. \ee

The matrix $k$-plane  is, in fact, a usual $km$-dimensional plane
in the Euclidean space $\Rnm$. To see this, we write $x=(x_{i,j})
\in\Ma$ and $t=(t_{i,j}) \in \Mt$ as column vectors
\be \bar x=\left(\begin{array} {c} x_{1,1}  \\ x_{1,2} \\ ... \\
x_{n,m}
\end{array} \right) \in \Rnm, \qquad \bar
t=\left(\begin{array} {c} t_{1,1}  \\ t_{1,2} \\ ... \\ t_{n-k,m}
\end{array} \right) \in \Rnkm,\ee
and denote \be \bar\xi=\diag (\xi, \ldots , \xi) \in
V_{nm,(n-k)m}.\ee Then (\ref{plane}) reads \be\label{4.5.1} \tau=
\tau(\bar \xi,\bar t)=\{\bar x\in\Rnm: \; \eqv\}. \ee The
$km$-dimensional planes (\ref{4.5.1})  form a subset of measure
zero in the affine Grassmann manifold of {\it all}
$km$-dimensional planes in $\Rnm$.

 The manifold $\Gr$ can be regarded as a fibre
bundle  the base of which is the
 ordinary Grassmann manifold $G_{n,k}$ of  $k$-dimensional linear subspaces $\eta$ of
$\rn$, and the fibres are homeomorphic to $\Mt$. Indeed, let
$\pi:\Gr\to G_{n,k} $ be the canonical projection which assigns to
each matrix $k$-plane $\tau(\xi,t)$ the subspace
\be\label{4.6}\eta=\eta(\xi)=\{y\in\rn:\xi ' y=0\} \in G_{n,k}.\ee
Let $\eta^\perp$ be the orthogonal complement of $\eta $ in $\rn$.
 The fiber $H_\eta=\pi^{-1}(\eta)$  is the set of
all matrix planes (\ref{plane}), when $t$ sweeps the space $\Mt$.

 Regarding $\Gr$  as a fibre bundle, one can utilize a
parameterization which is alternative to (\ref{plane}) and
one-to-one. Let \be x=[x_1 \dots x_m], \quad x_i\in\rn, \qquad
t=[t_1 \dots t_m], \quad t_i\in\rnk. \ee For $\tau= \tau(\xi,t)\in
\Gr$, we have
$$\tau=\{x\in\Ma:\xi ' x_i=t_i,\quad i=1\dots m \}.$$ Each
$k$-dimensional plane $\tau_i=\{x_i\in\rn:\xi ' x_i=t_i\}$ can be
 parameterized by the pair $(\eta, \lam_i)$, where $\eta$ is the subspace (\ref{4.6}), and
$\lam_i\in\eta^\bot$, $ \, i=1, \ldots , m,$ are columns of the
matrix $\lam=\xi t \in \Ma$. The corresponding
 parameterization
\be\label{p2} \tau =\tau (\eta,\lam), \qquad \eta \in G_{n,k},
\quad \lam=[\lam_1 \dots \lam_m], \quad \lam_i \in \eta^\perp, \ee
is one-to-one.

\subsection{Definition of the Radon transform}
The matrix $k$-plane Radon transform $\hat f$ of a function $f(x)$ on $\Ma$ assigns
to $f$ a collection of integrals of $f$ over all matrix planes $ \tau \in \Gr$.
 Namely, \[ \hat f (\tau)=\int_{x \in \tau} f(x).\] In order to
 give this integral precise meaning, we note that the matrix plane
$\tau= \tau(\xi,t), \; \xi\in\vnk$, $t\in\Mt$, consists of
``points"
\[ x=g_\xi \left[\begin{array} {c} u \\ t \end{array}
\right],  \] where $u \in \Mkm,$ and $ g_\xi \in SO(n)$ is a
rotation satisfying \be\label{4.24}
g_\xi\xi_0=\xi, \qquad \xi_0=\left[\begin{array} {c}  0 \\
I_{n-k} \end{array} \right] \in \vnk. \ee
 This observation leads to the following
\begin{definition}\label{rdef}
The Radon transform of a function $f(x)$ on $\Ma$ is defined as a
function  on the "matrix cylinder" $\cd$ by the formula
\be\label{4.9} \hat f (\tau) \equiv \rf=\intl_{\Mkm} f\left(g_\xi
\left[\begin{array} {c} u \\t
\end{array} \right]\right)du.
\ee
\end{definition}

 The reader is encouraged to check that  (\ref{4.9}) is independent of the choice of the
rotation $g_\xi:\xi_0\rightarrow\xi$. In terms of the one-to-one
parameterization (\ref{p2}), where $\tau =\tau (\eta,\lam)$, $\eta
\in G_{n,k}, \quad \lam=[\lam_1 \dots \lam_m] \in \Ma$, and $
\lam_i \in \eta^\perp$, the Radon transform (\ref{4.9}) reads \be
\hat f(\tau)=\intl_\eta dy_1\dots \intl_\eta f([y_1+\lam_1 \dots
y_m+\lam_m ]) \, dy_m. \ee If $m=1$, then $\rf$ is the ordinary
$k$-plane Radon transform that assigns to a function $f(x)$ on $
\bbr^n$ a collection of integrals of $f$ over all $k$-dimensional
planes [Hel]. A different definition of the matrix Radon
transform  was given by E.E. Petrov [Pe1]-[Pe3] (the case
$n-k=m$), and L.P. Shibasov [Sh1], [Sh2]  (the general case).

 The following
properties   can be easily checked.
\begin{lemma}\label{l4.3} Suppose that the Radon transform \[ f
(x) \longrightarrow \rf, \qquad x \in \Ma, \quad (\xi, t) \in \cd,
\]
 exists (at least almost everywhere). Then
 \be\label{ev} \hat f(\xi\theta ',\theta t)=\rf,  \qquad \forall \theta\in O(n-k). \ee
 If $ \; g(x)=\gamma x\beta +y \; $ where $ \;
\gamma\in O(n), \quad  \b\in O(m), \quad  y\in\Ma$, then
\be\label{4.23} (f \circ g)^\wedge (\xi, t)= \hat
f(\gamma\xi,t\beta +\xi'\gamma 'y). \ee In particular, if
 $f_y(x)=f(x+y)$, then \be
\label{4.4}
 \hat f_y(\xi ,t)= \hat f(\xi , \xi 'y+t).  \ee

\end{lemma}

The equality (\ref{ev}) is a matrix analog of the ``evenness
property" of the classical Radon transform, cf.  [Hel], p. 3.

\begin{lemma}\label{l4.2}\hskip 10truecm

\noindent {\rm (i)} If $f\in L^1(\Ma)$, then the Radon transform
$\rf$ exists for all $\xi\in\vnk$ and almost all $t\in\Mt$.
Furthermore,  \be \intl_{\Mt}\rf dt=\intl_{\Ma} f(x) \, dx, \qquad
\forall  \, \xi\in\vnk.\ee

\noindent {\rm (ii)} Let $||x||= (\tr(x'x))^{1/2}=(x_{1,1}^2+
\ldots +x_{n,m}^2)^{1/2}$. If $f$ is a continuous function
satisfying \be\label{4.21}
 f(x)=O(||x||^{-a}), \qquad a>km, \ee then $\rf$ exists  for all
 $\xi\in\vnk$ and all $t\in\Mt$.
\end{lemma}
\begin{proof} {\rm (i)} is a consequence of the Fubini
theorem: \bea \intl_{\Mt}\rf dt&=&\intl_{\Mt} dt \intl_{\Mkm}
f\left(g_\xi \left[\begin{array} {c} u \\t
\end{array} \right]\right)du \nonumber \\
 &=&\intl_{\Ma} f(g_\xi x) \,
dx=\intl_{\Ma} f(x) \, dx. \nonumber \eea (ii)
 becomes obvious if we regard $\tau=\tau (\xi, t)$  as a  $km$-dimensional plane  (\ref{4.5.1}) in
 $\Rnm$.
\end{proof}

A much deeper result is contained in the following
\begin{theorem}\label{t5.1}
If $f\in\lp$ then the Radon transform $\rf$  is finite for almost
all $(\xi, t) \in \cd$  provided \be\label{lp} 1\leq
p<p_0=\frac{n+m-1}{k+m-1}. \ee If $f$ is a continuous function
satisfying $ f(x)=O(|I_m+x'x|^{-\lam/2})$,  $ \lam>k+m-1$, then
$\rf$  is finite for  all $(\xi, t) \in \cd$.
\end{theorem}

 The proof of this theorem was given in [OR] using Abel type representation
 of the Radon transform of radial functions. The
conditions for $p$ and $\lam$ are sharp. For instance, one can
show [OR] that
\be\label{5.16}f_0(x)=|2I_m+x'x|^{-(n+m-1)/2p}(\log|2I_m+x'x|)^{-1}\ee
 belongs to $\lp$, and $\hat
f_0(\xi,t)\equiv\infty$ if  $p \ge p_0$. For $m=1$, the result of
Theorem \ref{t5.1} is due to Solmon [So]; see also [Ru2] for
another proof.

\subsection{Connection with the Fourier transform}
The Fourier transform of a function $f\in L^1(\Ma)$ is defined by
(\ref{ft}). The following statement is a matrix generalization of
the so-called Central Slice Theorem. It links together the
Fourier transform (\ref{ft}) and the Radon transform (\ref{4.9}).
In the case $m=1$, this theorem can be found in [Na, p. 11] (for
$k=n-1$) and [Ke, p. 283] (for any $0<k<n$).

For $y=[y_1 \dots y_m]\in\Ma$, let $\L (y)=\lin(y_1, \dots, y_m)$
be the linear hull of the $n$-vectors $y_1 \dots y_m$, that is the
smallest linear subspace containing $y_1 \dots y_m$. Suppose that
$\rank (y)=\ell$. Then $\dim\L (y)=\ell\leq m$.

 \begin{theorem}\label{CST}
 Let $f\in L^1(\Ma), \;n-k\geq m$. If  $y \in\Ma$, and
 $\zeta$ is an $(n-k)$-dimensional plane in $\bbr^n$ containing  $\L (y)$, then for
 any orthonormal frame $\xi\in\vnk$ spanning
 $\zeta$, there exists $b\in\Mt$ so that $y=\xi b$. In this case
 \be\label{4.19}
(\F f)(y)=\intl_{\Mt} \exp (i\, \tr(b't)) \, \rf \, dt, \ee
 or
\be\label{4.20} (\F f)(\xi b)=\F[\hat f(\xi,\cdot)](b), \quad
\xi\in\vnk, \quad b\in\Mt. \ee
 \end{theorem}

\begin{proof} Since each vector $y_j \; (j=1, \ldots , m) \; $ lies in $\zeta$, it
decomposes as   $y_j=\xi b_j$ for some $b_j \in \bbr^{n-k}$. Hence
  $y=\xi b$ where $b=[b_1 \ldots b_m]\in\Mt$. Thus it remains to  prove (\ref{4.20}). By (\ref{4.9}),
$$
\F [\hat f(\xi,\cdot)]( b)=\intl_{\Mt} \exp (i\, \tr(b't))  \, dt
\intl_{\Mkm} f\left(g_\xi \left[\begin{array} {c} u \\t
\end{array} \right]\right)du.$$
If $x=g_\xi \left[\begin{array} {c} u \\t
\end{array} \right]$, then, by (\ref{4.24}), $$\xi' x=\xi'_0 g'_\xi g_\xi \left[\begin{array} {c} u \\t
\end{array} \right]=\xi'_0 \left[\begin{array} {c} u \\t

\end{array} \right]=t, \qquad \xi_0= \! \left[\begin{array} {c}  0 \\
I_{n-k}
\end{array} \right] \! \in \vnk,
$$
and the Fubini theorem yields
$$
\F [\hat f(\xi,\cdot)](b)=\intl_{\Ma} \exp (i \, \tr(b'\xi 'x)) \,
f(x) \, dx=(\F f)(\xi b).
$$
\end{proof}
\begin{remark}\label{rem} It is clear that matrices $\xi$ and $b$ in
(\ref{4.19}) are not uniquely defined. In the case $\rank(y)=m$,
one can choose  $\xi$ and $b$ as follows. By taking into account
that $n-k \ge m$, we set
$$
u_0= \! \left[\begin{array} {c}  0 \\
I_{m}
\end{array} \right]  \! \in V_{n-k,m}, \qquad v_0= \! \left[\begin{array} {c}  0 \\
I_{m}\end{array} \right] \!  \in \vnm,
$$
so that $\xi_0 u_0=v_0$. Consider the polar decomposition $$y=v
r^{1/2},\qquad v\in\vnk, \quad r=y'y\in\p,$$ and let $g_v$ be a
rotation with the property $g_v v_0=v$. Then $$y=v r^{1/2}=g_v v_0
r^{1/2}=g_v \xi_0 u_0 r^{1/2}=\xi b,$$ where \be\label{4.25}
\xi=g_v\xi_0\in\vnk, \qquad  b=u_0 r^{1/2}\in\Mt.\ee
\end{remark}

\begin{theorem}\label{inj} {} \hfil

\noindent {\rm ($i$)} \ If $n-k\geq m$, then the Radon transform
$f\rightarrow\hat f$ is injective on the Schwartz space $\S(\Ma)$,
and $f$ can be recovered by the formula \bea \label {4.26} f(x)
&=& \frac{2^{-m}}{(2\pi)^{nm}}\intl_{\p} |r|^{\frac{n-m-1}{2}}dr
\nonumber \\ && \\ &\times& \intl_{\vnm} \exp (-i \,
\tr(x'vr^{1/2}))(\F\hat f(g_v \xi_0,\cdot))(\xi_0 'v_0 r^{1/2})dv.
\nonumber \eea

\noindent {\rm ($ii$)} \ For $n-k< m$, the Radon transform  is
non-injective.
\end{theorem}
\begin{proof}  By Theorem \ref{CST},  given the Radon transform  $\hat f$  of
 $f \in \S(\Ma)$,
the Fourier transform $(\F f)(y)$ can be evaluated  at each point
$y \in \Ma$ by the formula (\ref{4.19}), so that if $\hat f \equiv
0$ then $\F f \equiv 0$. Since $\F$ is injective, then $f \equiv
0$, and we are done.    Remark \ref{rem} allows to reconstruct $f$
from  $\hat f$,
 because (\ref{4.25})
 expresses $\xi$ and $b$ through $y \in \Ma$ explicitly. This gives (\ref{4.26}).
To prove {\rm ($ii$)}, we denote $$\frL_{n,m}=\{x\in\Ma:\rank
(x)=m\}.$$ This set is open in $\Ma$. Let $\psi$ be a
  Schwartz function with the Fourier transform
supported in $\frL_{n,m}$. By (\ref{4.20}), \be\label{4.14} \F
[\hat \psi(\xi,\cdot)](b)=\hat\psi(\xi b)=0 \quad \forall \;
\xi\in\vnk, \forall \; b\in\Mt, \ee because $\xi b \notin
\frL_{n,m}$ (note that since $n-k<m$, then $\rank(\xi b)<m$). By
injectivity of the Fourier transform in (\ref{4.14}),  we obtain $
\hat \psi(\xi,t)=0 \; \forall \xi, t$. Thus, for $n-k<m$, the
injectivity of the Radon transform fails.
\end{proof}

\begin{remark} After the paper had been  finished, we became aware of  another
 account of the topic of this subsection in [Sh1], [Sh2],
written in a different manner. Nevertheless, Remark \ref{rem},
formula (\ref{4.26}) and the statement (ii) of Theorem \ref{inj}
seem to be new. \end{remark}
\begin{remark} For $n-k>m$, the dimension of the manifold $G(n,k,m)$ of all
  matrix
$k$-planes in $\Ma$ is greater than that of the ambient space
$\Ma$, and the inversion problem is overdetermined. In the case
$n-k=m$ both dimensions coincide. The problem of reducing
overdeterminicity by fixing a certain ``invertible"
$mn$-dimensional complex of matrix planes, was studied in [Sh2].
\end{remark}

\section{The dual Radon transform and the Fuglede formula}

\setcounter{equation}{0}

\begin{definition} Let $\tau=\tau (\xi,t)$ be  a matrix
plane (\ref{plane}), $(\xi,t)\in\cd$. The  dual Radon transform
$\check \vp (x)$ assigns to a function $\vp(\tau)$ on $\Gr$ its
mean value over all matrix planes $\tau$ through $x$. Namely,
$$
\df=\intl_{\tau\ni x}\vp(\tau), \qquad x\in\Ma.
$$
This means that \bea\label{4.2} \qquad
\df&=&\frac{1}{\sigk}\intl_{\vnk}
\varphi(\xi,\xi'x)d\xi \\
&=&\intl_{SO(n)}\varphi(\g\xi_0,\xi_0'\g 'x)d\g, \qquad \xi_0=\left[\begin{array} {c}  0 \\
I_{n-k} \end{array} \right] \in \vnk. \nonumber \eea
\end{definition}

The mean value $\df$ apparently exists for all $x\in\Ma$ if
 $\vp$ is  a continuous function. Moreover,   $\df$
 is finite a.e. on $\Ma$ for any locally   integrable  function  $\vp$ [OR].
\begin{remark} The dual Radon transform $\df$ of a function  $\vp (\tau)
$, $\tau\in\Gr$, is independent of the parameterization $\tau=\tau
(\xi,t)$ in the sense that for any other parameterization
$\tau=\tau (\xi\theta',\theta t), \; \theta\in O(n-k),$  (see Sec.
\ref{s4.1}), (\ref{4.2}) gives the same result:
\[\frac{1}{\sigk}\intl_{\vnk}
\varphi(\xi\theta',\theta\xi'x)d\xi=\frac{1}{\sigk}\intl_{\vnk}
\varphi(\xi_1,\xi'_1 x)d\xi_1=\df.\]
 \end{remark}

\begin{lemma}
The duality relation \be\label{4.3} \intl_{\Ma} f(x)\df
dx=\frac{1}{\sigk} \intl_{\vnk}d\xi\intl_{\Mt}\fc\rf dt \ee holds
provided that either side of (\ref{4.3}) is finite for $f$ and
$\vp$ replaced by $|f|$ and $|\vp|$, respectively.
\end{lemma}
\begin{proof}
By (\ref{4.9}), the right hand side of (\ref{4.3})  is
\be\label{rs} \frac{1}{\sigk} \intl_{\vnk}d\xi\intl_{\Mt}\fc \, dt
\intl_{\Mkm} f\left(g_\xi \left[\begin{array} {c} u \\t
\end{array} \right]\right)du.\ee
Changing variables $x=g_\xi \left[\begin{array} {c} u \\t
\end{array} \right]$, we have
\[ \xi'x=(g_\xi \xi_0)'g_\xi \left[\begin{array} {c} u \\t
\end{array} \right]=\xi_0'\left[\begin{array} {c} u \\t
\end{array} \right]=t.\]
Hence,  by the Fubini theorem,  (\ref{rs}) reads \[
\frac{1}{\sigk} \intl_{\vnk}d\xi\intl_{\Ma} \varphi(\xi,\xi'x)
f(x) \, dx=\intl_{\Ma} \df  f(x) \, dx.\]
\end{proof}

Now we state the main theorem.
\begin{theorem}\label{main} Let $f\in\lp, \quad 1\leq
p<(n+m-1)/(k+m-1)$, or  $f$ is a continuous function
satisfying $ f(x)=O(|I_m+x'x|^{-\lam/2})$,  $ \lam>k+m-1$. Then
 \be\label{fu} c(\hat
f)^{\vee} (x) \! = \! (I^k f)(x), \qquad c \! = \! 2^{-km}
\,\pi^{-km/2} \, \Gam_m\Big ( \frac{n \! - \! k}{2}\Big ) /
\Gam_m\Big ( \frac{n}{2}\Big ), \ee (the generalized Fuglede
formula).
\end{theorem}
\begin{proof} The proof given below is applicable to any
locally integrable function $f$ for which either side of
(\ref{fu}) is finite provided $f$ is replaced by $|f|$. Let $f_x
(y)=f(x+y)$. By (\ref{4.2}) and (\ref{4.4}), \bea(\hat f)^{\vee}
(x)&=&\frac{1}{\sigk}\intl_{\vnk} \hat
f(\xi,\xi'x)d\xi=\frac{1}{\sigk}\intl_{\vnk} \hat f_x(\xi,0)d\xi
\nonumber \\
&=&\intl_{\frM_{k,m}} du \intl_{SO(n)} f\left ( x + \gam
\left[\begin{array} {c} u \\ o\end{array} \right]\right ) d\gam.
\nonumber \eea  This coincides with $c^{-1}(I^k f)(x)$ for $k<m$,
see Definition \ref{cor}. If $k \ge m, \; d=(m+1)/2$, we pass to
polar coordinates and get \bea  (\hat f)^{\vee}
(x)&=&2^{-m}\intl_{V_{k,m}} dv\intl_{\P_{m}}|r|^{k/2 -d} dr
\intl_{SO(n)} f_x\left ( \gam \left[\begin{array} {c} vr^{1/2} \\
0\end{array} \right]\right ) d\gam, \nonumber \\ &=& \frac{2^{-m}
\, \sig_{k,m} }{\sig_{n,m}}\intl_{\P_{m}}|r|^{k/2 -d} dr
\intl_{V_{n,m}} f_x (wr^{1/2})dw \nonumber \\ &=& \frac{
\sig_{k,m} }{\sig_{n,m}}\intl_{\Ma} f(x+y) |y'y|^{(k-n)/2}
dy\nonumber \\ &=&c^{-1}(I^k f)(x). \nonumber \eea \end{proof}

\begin{corollary}\label{kint} {\rm (cf. Conjecture 3.4)} For
$f \in L^p (\Ma)$ and  $k \in \bbn$, the  Riesz potential $(I^k f)(x)$
is finite a.e. on $\Ma$ if and only if \be \label{co} 1\leq
p<(n+m-1)/(k+m-1). \ee If $f$ is a continuous function satisfying
$ f(x)=O(|I_m+x'x|^{-\lam/2})$, then  $(I^k f)(x)$  is finite for
all $x \in\Ma$ if and only if $ \lam>k+m-1$.
\end{corollary}

This statement holds by Theorem \ref{t5.1}.

\section{Inversion problem for the Radon transform}

\setcounter{equation}{0}

The Fuglede formula $c(\hat f)^{\vee} =I^k f$ reduces the
inversion problem for the Radon transform to Riesz potentials.
This is exactly the same situation as in the rank-one case.
Whereas for the ordinary $k$-plane transform and Riesz potentials
a variety of pointwise inversion formulas is available in a large
scale of function spaces [Ru1], [Ru2], in the higher rank case we
cannot obtain pointwise  inversion formulas rather than  for
Schwartz functions via the Fourier transform  (see (4.22)) or in terms of
divergent integrals  understood somehow in the regularized sense.
This is still an open problem.

Below we show how the unknown ``rough" function $f$ can be
recovered in the framework of the theory of distributions. First
we specify the  space of test functions. From the Fourier
transform formula $(h_\a, f)=(|y|_m^{-\a},(\F f)(y))$, it is
evident that the Schwartz class $\S \equiv \S(\Ma)$ does not suit
well enough
  because it is not invariant under multiplication by
$|y|_m^{-\a}$. To get around this difficulty, we follow
 an  idea of V.I.
Semyanistyi [Se] suggested for  $m=1$. Let $\Psi \equiv \Psi(\Ma)$
be the subspace of functions $\psi (y) \in \S$ vanishing on the
set \be\label{sets} \{y: \, y \in \Ma, \;  \rank (y)<m \}=\{y: \,y
\in \Ma, \; |y'y|=0 \}\ee with all derivatives (the coincidence of
both sets in (\ref{sets}) is clear because $\rank (y)=\rank
(y'y)$, see, e.g., [FZ], p. 5). The set $\Psi$ is a closed linear
subspace of $\S$. Therefore, it can be regarded as a linear
topological space with the induced topology of $\S$. Let $
\Phi\equiv\Phi (\Ma)$ be the
 Fourier image of $\Psi$. Since the Fourier transform
$\F$ is an automorphism of $\S$ (i.e., a topological isomorphism
of $\S$ onto itself), then $\Phi$ is a closed linear subspace of
$\S$. Having been equipped with the induced topology of $\S$, the
space $\Phi$ becomes a linear topological space isomorphic to
$\Psi$ under the Fourier transform. We denote by  $
\Phi'\equiv\Phi' (\Ma)$ the space of all linear continuous
functionals (generalized functions) on  $\Phi$.  Since for any
complex $\a$, multiplication by $|y|_m^{-\a}$ is an automorphism
of $\Psi$, then, according to the general theory [GSh], $I^\a$, as
a convolution with $h_\a$, is an automorphism of $\Phi$, and we
have $$\F[I^\a f](y)=|y|_m^{-\a}\F[f](y)$$ for all
$\Phi'$-distributions $f$.

In the  rank-one case, the spaces $\Phi$, $\Psi$, their duals and
generalizations were  studied by P.I. Lizorkin, S.G. Samko and
others in view of applications to the theory of function spaces
and fractional calculus; see [Sa], [SKM], [Ru1] and references
therein.

The Fuglede formula (\ref{fu}), Theorem \ref{t5.1}, and Corollary
\ref{kint} imply the following
\begin{theorem} Let $f \in L^p (\Ma), \; 1\leq p<(n+m-1)/(k+m-1)$
or $f$ is a continuous function satisfying $
f(x)=O(|I_m+x'x|^{-\lam/2})$ for some $ \lam>k+m-1$. Then the
Radon transform $g=\hat f$ is well defined, and $f$ can be
recovered from $g$ in the sense of $\Phi'$-distributions by the
formula \be (f,\phi)=c(\check g, I^{-k}\phi), \qquad \phi \in
\Phi, \ee where $$(I^{-k}\phi)(x)=(\F^{-1}|y|_m^{k}\F\phi)(x),
\quad c=2^{-km} \,\pi^{-km/2} \, \Gam_m\Big ( \frac{n-k}{2}\Big )/
 \Gam_m\Big ( \frac{n}{2}\Big ).$$
\end{theorem}
\begin{remark}\label{clap} For $k$ even, the Riesz potential $I^kf$ can be
inverted (in the sense of $\Phi'$-distributions) by repeated
application  of the Cayley-Laplace operator $ \Del_m =\det
(\partial' \partial), \;  \partial=(\partial/\partial x_{i,j})$ [Kh]. In the Fourier
terms, this operator agrees with multiplication by $(-1)^m
|y|_m^2$, and therefore, $(-1)^m \Del_m I^\a f=I^{\a-2} f$ in the
$\Phi'$-sense.
\end{remark}

\end{document}